\documentclass[reqno,11pt]{amsart}

\newtheorem{theorem}{Theorem}[section]
\newtheorem{lemma}[theorem]{Lemma}
\newtheorem{corollary}[theorem]{Corollary}

\theoremstyle{definition}

\theoremstyle{remark}
\newtheorem{remark}[theorem]{Remark}

\newcommand{\mysection}[1]{\section{#1}
\setcounter{equation}{0}}

\newcommand{\bR}{\mathbb R}

\newcommand\cL{\mathcal{L}}

\begin{document}
\title[Smooth functions with values in polyhedra] {On
factorizations of smooth nonnegative matrix-values
functions and on
smooth functions with values in polyhedra}

\author[N.V. Krylov]{N.V. Krylov}

\address
{127 Vincent Hall, University of Minnesota, Minneapolis, MN
55455, USA}
\thanks{The work  was partially supported by
NSF Grant DMS-0653121}
\email{krylov@math.umn.edu}



\begin{abstract}
We discuss the possibility to represent smooth
nonnegative matrix-valued functions as  finite linear
combinations
of fixed matrices with positive real-valued coefficients
whose square roots are Lipschitz continuous.
This issue is reduced to a similar problem
for smooth functions with values in a polyhedron.
\end{abstract}

\maketitle

\mysection{Motivation}

One of the main goals of the article is to understand
what kind of optimal control problems
of diffusion processes is covered by the results
of \cite{DK} and \cite{Kr}, where the processes are given by
It\^o equations in a ``special" form, such that in the
corresponding
Bellman equation the second order part is represented as 
the sum of second-order derivatives
with respect to fixed vectors (independent of the control
parameter) times squares of real-valued functions
that are   Lipschitz continuous with respect to the
space variables. Roughly speaking the answer is that
all control problems with twice continuously differentiable
diffusion matrices fall into the scheme of \cite{DK} and \cite{Kr}
whenever property (A) holds: these matrices for all values of
control and time and space variables belong to a fixed
polyhedron in the set of symmetric nonnegative matrices.
In the author's opinion the control problems
with property (A) are the only ones which
admit finite-difference approximations with
monotone schemes based on scaling of a fixed mesh.

For functions $w(z)$ given in a Euclidean space and vectors
$\xi$ in that space set
$$
w_{(\xi)}=(\xi,\nabla w)=\sum_{i}\xi^{i}w_{z^{i}},\quad
w_{(\xi)(\xi)}=\sum_{i,j}\xi^{i}\xi^{j}w_{z^{i}z^{j}} .
$$

In many situations one needs to represent a 
$d \times d  $ nonnegative symmetric
 matrix $u$ as the square of a matrix 
or more generally as the product $v v^{*} $, where $v$
is not necessarily a square matrix. 
If $u=(u^{ij})=vv^{*}$ and $v=(v^{ik})$ and 
for each $k$ we introduce the
vector  $v^{k}=(v^{ik})\in\bR^{d}$, then
for any smooth $f(x)$ given on $\bR^{d}$ 
and the operator
\begin{equation}
                                            \label{5.31.01}
Lf:=\sum_{i,j}u^{ij}f_{x^{i}x^{j}}
\end{equation}
we have
\begin{equation}
                                            \label{5.18.1}
Lf=\sum_{k}f_{(v^{k})(v^{k})}.
\end{equation}
In fact, as is easy to see having \eqref{5.18.1}
for all smooth $f$ is equivalent
to the validity of the formula $u=vv^{*}$.

There are very many $v$ such that $u =vv^{*}$
 and then a few questions arise:

(i) if $u$ is a measurable function of a parameter,
 can one find a measurable
$v$?

(ii) if $u$ is smooth, can one find a Lipschitz continuous
$v$?

The answer to the first question is easy and positive.
Indeed, one can take
$$
v=c\int_{0}^{\infty}t^{-3/2}(e^{-u t}-1)\,dt,
$$
where $c$ is an appropriate constant. This defines $v$
as the square root of $u$. Since long ago it is known that
the square root
of a twice differentiable nonnegative matrix-valued function  is
 Lipschitz continuous (see \cite{Fr}, \cite{PS}).
This result was used in the investigation
of  solvability of degenerate elliptic and parabolic
second-order equations by using probabilistic or
classical approaches.

However, there are applications in which 
formula \eqref{5.18.1} is not very convenient.
One of these applications is related to
finite-difference approximations of solutions
to elliptic and parabolic equations with variable
coefficients $u^{ij}$.
Formula \eqref{5.18.1} suggests replacing $f_{(v^{k})(v^{k})}$
with the second-order difference along vectors
$v^{k}$ and if $v^{k}$ vary, it may be impossible to find
a reasonable mesh on which the approximation
operator would make sense. 
This problem does not arise if $v^{k}=\sigma_{k}\gamma^{k}$,
where $\gamma^{k}$ are constant vectors and $\sigma_{k}$
are real-valued functions, because then
\begin{equation}
                                            \label{5.18.2}
Lf =\sum_{k}
\sigma_{k}^{2}f_{(\gamma^{k})(\gamma^{k})}
\end{equation}
and one can concentrate on meshes that are obtained
by contracting 
$$\{\sum_{k}n_{k}\gamma_{k}:n_{k}=0,\pm1,...\}.
$$
According to Remark 2.1 of \cite{DK}
considering operators $L$ in   form \eqref{5.18.2}
is rather realistic
from the point of view of numerical approximations.
It turns out
that if we fix a finite subset $B\subset\bR^{d}$,
such that $\text{Span}\,B=\bR^{d}$,
      and if $L$ from \eqref{5.31.01} admits a finite-difference
approximation
$$
L_{h}f(0)=\sum_{y\in B}p_{h}(y)f(x+hy)\to Lf(0)
\quad\text{as}\quad h\downarrow0,
\quad\forall f\in C^{2}
$$
and $L_{h}$ are monotone, that is $L_{h}f(0)\geq0$
whenever $f(x)\geq f(0)$ on $\bR^{d}$,
 then automatically
    $L$ is written in the form \eqref{5.18.2}
with some $\gamma^{k}\in B$.

Now the question is: If $u=u(x)$, under which  assumptions
can one find constant vectors $\gamma_{k}$'s and functions
$\sigma_{k}^{2}(x)$ in order for \eqref{5.18.2} to hold?
Perhaps, Motzkin and Wasow (see \cite{MW})
were the first to address this question
in the framework of finite-difference approximation.
They proved (see also Lemma 17.13 of \cite{GT})
that if we denote by $S[\lambda,\Lambda]$ the (closed) set of
positive $d \times d $ matrices with eigenvalues
lying in the interval $[\lambda,\Lambda]$, where
$0<\lambda\leq\Lambda$, then there exist  a finite set
of unit vectors $\gamma^{1},...,\gamma^{m}\in\bR^{d_{2}}$
and  numbers $0<\lambda^{*}<\Lambda^{*}$, such that
any $u\in S[\lambda,\Lambda]$ can be written in the form
\begin{equation}
                                            \label{5.18.3}
u^{ij}=\sum_{k=1}^{m} \beta_{k}\gamma ^{ik}\gamma ^{jk},
\end{equation}
where the numbers $\beta_{k}$ satisfy the inequalities
$\lambda^{*}\leq\beta_{k}\leq\Lambda^{*}$. In that case
\eqref{5.18.2} holds with $\sigma^{2}_{k}=\beta_{k}$.
This fact was used in the development of the theory
of {\em fully\/} nonlinear second-order elliptic and parabolic
equations.

One can give a quite easy explanation of this result.
If we   take any $\lambda_{1}<\lambda$ and $\Lambda_{1}>\Lambda$,
the set $S[\lambda_{1},\Lambda_{1}]$ will contain
an open polyhedron $P^{o}$ containing $S[\lambda ,\Lambda ]$.
Each point of a polyhedron is represented as a convex
combination of its vertices and one easily obtains
\eqref{5.18.3}, for instance, as in the proof of Lemma 5.5.4
of \cite{Kr1}.

With a little more effort one can get more convenient
representations.
We claim that given any open bounded
 polyhedron $P^{o}$ in a Euclidean
space $E$ of points $y$
with vertices, say $y_{1},....,y_{n}$, there
exist
 infinitely differentiable  
functions $p_{k}(y)>0$, $k=1,...,n$, such that
for any $y\in P^{o}$
\begin{equation}
                                            \label{5.18.4}
y=\sum_{k}p_{k}(y)y_{k},\quad\sum_{k}p_{k}(y)=1.
\end{equation}
This fact is proved by induction on the dimension of $P^{o}$.
First, without losing generality one may assume that
the volume of $P^{o}$ in $E$ is strictly greater than zero.
Then, assume that the fact is true for any face
of $P^{o}$ and then solve Laplace's equation $\Delta p_{k}=0$
in $P^{o}$ with boundary condition $p_{k}(y)=\bar{p}_{k}(y)$
on   $\partial P^{o}$, where $\bar{p}_{k}$
is the weight of the vertex $y_{k}$ in the representation
of $y\in\partial P^{o}$, which is supposed to hold by the
induction hypothesis. 
Of course, if $y\in\partial P^{o}$ and $y_{k}$ do not belong to the
same face, we set $\bar{p}_{k}(y)=0$.
 Then by the well-known properties of
harmonic functions $p_{k}>0$ in $P^{o}$, they are infinitely
differentiable
in $P^{o}$ and since 
$$
\Delta \sum_{k}p_{k}(y)y_{k}=0,\quad\text{in}\quad P^{o}
\quad\text{and}
\quad 
\sum_{k}p_{k}(y)y_{k}=y\quad\text{on}\quad\partial P^{o}
$$
and $\Delta y=0$, by uniqueness
we have the first relation in \eqref{5.18.4}. The second one is
obtained similarly from the fact that it holds on $\partial P^{o}$
and $\Delta 1=0$.

After having proved the claim we return to the original $P^{o}$
and write for any $u\in P^{o}$
\begin{equation}
                                            \label{5.19.1}
u=\sum_{k}p_{k}(u)u_{k},\quad
\sum_{k}p_{k}(u)=1,
\end{equation}
where $u_{k}\in P^{o}\subset S[\lambda_{1},\Lambda_{1}]$,
$p_{k} $  are infinitely differentiable in $P^{o}$,
in particular, in $S[\lambda,\Lambda]$,
$p_{k}>0$ in $P^{o}$, in particular, they are bounded away from zero
on the closed set $S[\lambda,\Lambda]$.
Now to obtain \eqref{5.18.3} from \eqref{5.19.1}
it only remains to recall that if $\xi_{ki}$, $i=1,...,d$,
are unit eigenvectors of $u_{k}$ with eigenvalues $\mu_{ki}$,
then $\lambda_{1}\leq\mu_{ki}\leq\Lambda_{1}$ and
$$
u_{k}=\sum_{i}\mu_{ki}\xi_{ki}\xi^{*}_{ki}.
$$

The above construction of $p_{k}(u)$ has a very substantial advantage
over the original one in \cite{MW} (or \cite{GT}
and \cite{Kr1}). Namely, it is seen
that if $u=u(x)$ is a smooth function of a parameter $x$, then in the
representation
\begin{equation}
                                            \label{5.19.2}
u(x)=\sum_{k}p_{k}(u(x))u_{k}
\end{equation}
or in the implied representation 
\eqref{5.18.2}
the functions $p_{k}(u(x))$, $p^{1/2}_{k}(u(x))$,
and $\sigma _{k}(x)$ are as smooth as $u(x)$ is.

We see that from the point of view of the possibility
of applying
numerical approximations
 to uniformly nondegenerate equations the situation 
looks quite promising.
For degenerate equations and fully nonlinear
equations the situation is much more
complex. In this case we again may try to prove
\eqref{5.19.1} with $p_{k}$ such that $p_{k}^{1/2}(u)$
is Lipschitz continuous in $u$. However, this
is impossible even if $d=1$ and $S[\lambda,\Lambda]
=[0,1]$. In this case, naturally $u_{1}=0$, $u_{2}=1$, and
$p_{2}(u)=u$, so that $p_{2}^{1/2}(u)$ is not
Lipschitz continuous.

On the other hand, in numerical approximation
or probabilistic approach one needs 
$p^{1/2}_{k}(u(x))$ to be Lipschitz continuous function of
$x$, rather than $p^{1/2}_{k}(u )$ to
be Lipschitz continuous function of
$u$. This slight difference makes the problem solvable
in some cases. For instance, in the above case that $d=1$
it is known that for any nonnegative
{\em twice\/} continuously differentiable function
$u(x)$ its square root $u^{1/2}(x)$ is
 Lipschitz continuous.

Another example is given by the functions with values
in the set of the so-called diagonally dominant
nonnegative symmetric matrices, which are quite popular
in the literature (see \cite{BJ}, \cite{KD}).
  These are the ones with the property
\begin{equation}
                                            \label{5.19.3}
2u^{ii}\geq\sum_{j=1}^{d}|u^{ij}|,\quad  i=1,...,d.
\end{equation}
Let $D$ be the set of symmetric matrices satisfying \eqref{5.19.3}
and such that $\text{trace}\,u=1$. 
The author heard some doubts that, say the results of
\cite{DK} are applicable to equations whose variable
coefficients
of second order derivatives form matrices of class $D$.
The point is that the equations in \cite{DK} are assumed to have
the structure associated with \eqref{5.18.2}
with Lipschitz continuous $\sigma_{k}$. A naive way fails
to take $e_{1},...,e_{d}$ as the
standard basis vectors in $\bR^{d}$
and
write a representation like \eqref{5.18.3} for a $  D$-valued
function $u(x)$ as
(see \cite{BJ}) 
$$
u(x)=\sum_{i\ne j}\big[( u^{ii}(x)-|u^{ij}(x)|)e_{i}e^{*}_{i}
+(1/2)(u^{ij}(x))^{+}(e_{i}+e_{j})(e_{i}+e_{j})^{*}
$$
$$
+(1/2)(u^{ij}(x))^{-}(e_{i}-e_{j})(e_{i}-e_{j})^{*}\big],
$$
where we used the notation $a^{\pm}=(1/2)(|a|\pm a)$.
The reason for the failure is that no smoothness
assumptions on $u(x)$ can guarantee that
$[(u^{ij})^{+}]^{1/2}$ is Lipschitz continuous
for $i\ne j$. One needs a nontrivial structural assumption
for that.

Nevertheless, in \cite{Kr} for $d=2$ the author 
gave explicit formulas for representing 
twice continuously differentiable $D$-valued functions 
in the form \eqref{5.19.2} with Lipschitz continuous
$p_{k}^{1/2}(u(x))$. The result of the present article shows that
such representation exists for any $d$.
In addition, it turns out that the set of 
diagonally dominant matrices
can be replaced with any set which is a polyhedron
in the set of  $d\times d$ matrices.
 By the way, observe that \eqref{5.19.3}
can be equivalently written as
$$
2u^{ii}\geq\sum_{j=1}^{d} \varepsilon^{ij}u^{ij} ,\quad  i=1,...,d,
\varepsilon^{ij}=\pm1.
$$
In $D$ we also have $\text{trace}\,u=1$.
Therefore, the bounded set $D$ is described
by means of finitely  many linear equalities and inequalities,
and hence $D$ is a polyhedron in the space of $d\times d$ matrices.
Speaking about the case that $d=2$, it is also 
worth noting that in \cite{BOZ} an efficient
algorithm is introduced for approximating arbitrary
$2\times 2$ nonnegative matrices with matrices of the form
$\sum_{k}p_{k}\xi_{k}\xi_{k}^{*}$, where $\xi_{k}\in \bR^{2}$.

Finally, we reiterate that representation \eqref{5.19.2}
leads to \eqref{5.18.2} and the latter means that
we have the following factorization:
 $$
u=vv^{*},\quad\text{where}\quad
 v^{ik}=\sigma_{k}\gamma^{ik}.
$$ 

Starting from this point we forget about matrices
and work with functions having values in a polyhedron.
Our main results are presented in Section \ref{section 5.21.1},
Theorem \ref{theorem 11.10.2} is proved in
Section \ref{section 5.21.1} and Theorem \ref{theorem 07.5.7.1}
is proved in Section \ref{section 5.21.3}. Section 
\ref{section 5.21.2} contains an investigation of an
auxiliary function some additional
information on which is provided
in Section \ref{section 5.21.4}.

The author discussed the article with  Hongjie Dong whose comments
are greatly appreciated.

\mysection{Main results}
                                            \label{section 5.21.1}

 Let $P$ be a closed bounded convex polyhedron in $\bR^{d}$ with 
distinct vertices
$a_{1},...,a_{n}$, where $n\geq2$. Let $d_{1}\geq1$
be an integer.

\begin{theorem}
                                              \label{theorem 11.10.2}
Let $u(y)$ be a  $P$-valued functions defined on $\bR^{d_{1}}$.
Assume that the first and second order derivatives of $u$
are bounded and continuous on $\bR^{d_{1}}$.
Then there exist real valued nonnegative functions $u_{1}(y),
...,u_{n}(y)$ such that
\begin{equation}
                                               \label{11.10.5}
\sum_{k}u_{k}(y)\equiv1,\quad u(y)\equiv
\sum_{k}u_{k}(y)a_{k},
\end{equation}
and $u^{1/2}_{k}$ are Lipschitz continuous on $\bR^{d_{1}}$
with a constant which depends only on $P$ and 
$\sup\{|u_{(\eta)(\eta)}(y)|:|\eta|=1,y\in\bR^{d_{1}}\}$.
\end{theorem}

Clearly the following assumption which we keep throughout
the paper does not restrict generality:
$$
a_{1}=0,\quad\text{Span}\,(a_{2},...,a_{n})=\bR^{d}.
$$

To prove Theorem \ref{theorem 11.10.2} we need the following result.
For $\xi\in\bR^{d}\setminus\{0\}$ and $x\in P$
denote by $d(x,\xi)$ the distance from $x$ to $\partial P$
along the ray $x+t\xi/|\xi|$, $t\geq0$. Introduce,
$P^{o}$ as the interior of $P$.

Denote by $\Phi$ the set of $d-1$-dimensional faces of $P$ and for
$\Gamma\in\Phi$ and $x\in P^{o}$
introduce $d_{\Gamma}(x)$ as the distance from  $x$
to $\Gamma$. Also let $n_{\Gamma}$ be a unit normal
vector to $\Gamma$.

\begin{theorem}
                                       \label{theorem 07.5.7.1}

On $P$ there exist Lipschitz continuous nonnegative
functions $p_{1}(x),...,p_{n}(x)$ which are infinitely
differentiable in $P^{o}$ and such that

(i) $p_{k}>0$ in $P^{o}$;

(ii) in $P$ we have
$$
\sum_{k}p_{k}(x)=1,\quad x=\sum_{k}p_{k}(x)a_{k};
$$

(iii) for any  $\xi\in\bR^{d} $ and $x\in P^{o}$ we have
\begin{equation}
                                       \label{07.5.7.1}
\frac{|p_{k(\xi)}(x)|}{
p_{k}^{1/2}(x)}\leq N\max_{\Gamma\in\Phi}
\frac{|(n_{\Gamma} ,\xi)|}{d^{1/2}_{\Gamma}(x)}\,,
\end{equation}
where $N$ is a finite constant depending only on $P$.
\end{theorem}

{\bf Proof of Theorem \ref{theorem 11.10.2}}. Take a point
$x_{0}\in P^{o}$ and for 
 $t\in(0,1)$ set $u_{t}(y)=tu(y)+(1-t)x_{0}$.
Then $u_{t}$ takes values in $P^{o}$. Assume that
for $P^{o}$-valued functions the statement of
Theorem \ref{theorem 11.10.2} is true. Then,
for each $t\in(0,1)$ 
there exist real valued nonnegative functions $u_{t1}(y),
...,u_{tn}(y)$ such that
$$
\sum_{k}u_{tk}(y)\equiv1,\quad u_{t}(y)\equiv
\sum_{k}u_{tk}(y)a_{k},
$$
and $u^{1/2}_{tk}$ are Lipschitz continuous on $\bR^{d_{1}}$
with a constant independent of $t$. By the Arzel\`a-Ascoli theorem
it follows that there exists a sequence $t_{n}\uparrow1$ such that
$u_{t_{n}k}(y)$ converge to some functions $u_{k}(y)$ for
each $y$ and $u_{k}^{1/2}$ are Lipschitz continuous.
Obviously, these are the functions which we need.
We see that without losing generality we may assume that
$u(y)\in P^{o}$ for all $y$.

We will be using the well-known fact that if
we have a nonnegative twice continuously differentiable
function $f(y)$ given on $\bR^{d_{1}}$ and having
bounded second-order derivatives, then for any $y\in\bR^{d_{1}}$
$$
|\nabla f(y)|^{2}\leq 4f (y)\sup\{|f_{(\eta)(\eta)}(z):
|\eta|=1,z\in\bR^{d_{1}}\} .
$$
 
Now
take $p_{k}$ from Theorem \ref{theorem 07.5.7.1} and set
$u_{k}(y)=p_{k}(u(y))$. Then 
the equations \eqref{11.10.5} obviously hold.
Since $p_{k}$ are positive and
infinitely differentiable in $P^{o}$, $u_{k}$ are positive and  
continuously differentiable in $\bR^{d_{1}}$. Therefore,
to estimate the Lipschitz constant of $u^{1/2}_{k}$
it suffices to estimate its first order directional derivatives.

Fix a $y,\eta\in \bR^{d_{1}}$ with $|\eta|=1$  and set $x=u(y)$,
$\xi=u_{(\eta)}(y)$. Then by \eqref{07.5.7.1}
\begin{equation}
                                                  \label{11.10.7}
2|(u^{1/2}_{k})_{(\eta)}(y)|=\frac{|p_{k(\xi)}(x)|}
{p^{1/2}_{k}(x)}\leq N\max_{\Gamma\in\Phi}
\frac{|(n_{\Gamma} ,\xi)|}{d^{1/2}_{\Gamma}(x)}\,.
\end{equation}

Next, take a face $\Gamma\in\Phi$ and let it be given as
$\{x:(n_{\Gamma},x)=b\}$, where $b$ is a constant. By multiplying
  $n_{\Gamma}$ and $b$ by $-1$ if needed we may assume that
$$
(n_{\Gamma},w)\geq b\quad\forall w\in P.
$$
Then $f(z):=(n_{\Gamma},u(z))-b$ is a nonnegative
twice continuously differentiable function
on $\bR^{d_{1}}$. By the above
$$
|(n_{\Gamma},\xi)|^{2}=|f_{(\eta)}(y)|^{2}\leq
Nf(y)=N|(n_{\Gamma},x)-b|=Nd_{\Gamma}(x),
$$
where
$$
N=4\sup_{|\eta|=1,z}|u_{(\eta)(\eta)}(z)|.
$$
This and \eqref{11.10.7} bring the proof of the theorem to
an end.
 
\mysection{An auxiliary function}
                                            \label{section 5.21.2}
 For $x\in P^{o}$ define
$$
U(x)=\max\{\sum_{i=1}^{n}\ln p_{i}: p_{i}>0,
\sum_{i=1}^{n}p_{i} =1,\quad\sum_{i=1}^{n}p_{i} a_{i}=x\}.
$$
Obviously, $U\leq0$
and for each $x\in P^{o}$ there exists
$p_{1},...,p_{n}$ achieving the maximum.
 
\begin{lemma}  
(i) The function $U$ is strictly concave and therefore continuous
  in   $P^{o}$.

(ii) For each $x\in P^{o}$ there exists a 
unique set $p_{1}(x),...,p_{n}(x)>0$
such that
$$
\sum_{i=1}^{n}p_{i} (x)=1,\quad\sum_{i=1}^{n}p_{i} 
(x)a_{i}=x,\quad U(x)=\sum_{i}\ln p_{i}(x).
$$

(iii) The functions $p_{1}(x),...,p_{n}(x)$ are continuous in
$P^{o}$.

\end{lemma}

Proof. (i)
Take $x,y\in P^{o}$, $t,s\in(0,1)$, such that $t+s=1$,
and let
$p_{1},...,p_{n}$ and $q_{1},...,q_{n}$ be some sets achieving the
maximums for $x$ and $y$ respectively. Then for
$r_{i}=tp_{i}+sq_{i} $ we have
$$
\sum_{i}r_{i}=1,\quad\sum_{i}r_{i}a_{i}=tx+sy .
$$
Hence,
$$
U(tx+sy)\geq\sum_{i}\ln (tp_{i}+sq_{i}) 
\geq  t\sum_{i} \ln p_{i}+s\sum_{i}\ln q_{i} 
$$
$$
=tU(x)+sU(y) ,
$$
where the second inequality is strict if
$p_{i}\ne q_{i}$   for at least one $i$.
This is certainly the case if $x\ne y$, which proves (i).
Another case would appear if $x=y$ and we assumed that
there are two different sets $p_{1},...,p_{n}$
and $q_{1},...,q_{n}$ achieving $U(x)$. But then the above
computations would lead to a wrong conclusion that
$U(x)>U(x)$. This proves (ii).

Finally (iii) follows from the continuity of $U(x)$
and assertion (ii). The lemma is proved.

\begin{lemma}
                                               \label{lemma 5.15.1}
The function $U$ is continuously differentiable in $P^{o}$ and
for any $x\in P^{o}$ and $\xi\in\bR^{d}$, which is represented as
\begin{equation}
                                                \label{5.17.1}
\xi=\sum_{k}q_{k}(x-a_{k})
\end{equation}
with some numbers $q_{k}$, we have
\begin{equation}
                                                \label{8.23.1}
\sum_{k}\frac{q_{k}}{p_{k}(x)}=n\sum_{k}q_{k}-U_{(\xi)}(x).
\end{equation}
In particular, as $\xi=x-a_{k}$,
\begin{equation}
                                                \label{5.31.2}
\frac{1}{p_{k}(x)}=n- U_{(x-a_{k}) }(x).
\end{equation}
\end{lemma}

Proof. Fix an $x_{0}\in P^{o}$ and let $\lambda\in\bR^{d}$
be such that the graph of the function
$ (\lambda,x-x_{0})+U(x_{0})$ is a supporting plane for the graph of
$U(x)$ at $(x_{0},U(x_{0}))$.
Set $b:=\sum_{k}q_{k}$   and write
$$
x_{0}+t\xi=\sum_{k}[(1+bt)p_{k}(x_{0})- t
q_{k}]a_{k}.
$$
For sufficiently small $t$ we have
$(1+bt)p_{k}(x_{0})-t
q_{k}>0$ and
$$
\sum_{k}[(1+bt)p_{k}(x_{0})-t
q_{k}]=1.
$$
It follows that for small $t$
$$
t(\lambda,\xi) +U(x_{0})\geq
U(x_{0}+t\xi)\geq\sum_{k}\ln
[(1+bt)p_{k}(x_{0})-t
q_{k}]
$$
with equalities instead of the inequalities for $t=0$.
By differentiating at $t=0$ the extreme terms we find
$$
(\lambda,\xi) =n\sum_{k}q_{k}-\sum_{k}\frac{q_{k}}{p_{k}(x_{0})}.
$$
If there is another vector $\mu\in\bR^{d}$
such that the graph of the function
$ (\mu,x-x_{0})+U(x_{0})$ is a supporting plane for
the graph of
$U(x)$ at $(x_{0},U(x_{0}))$, then the above formula
implies that $\lambda-\mu\perp\xi$. This holds for 
any $\xi$ admitting representation
\eqref{5.17.1} with $x_{0}$ in place of $x$. Since
\begin{equation}
                                                       \label{5.21.1}
\text{Span}\,\{(x_{0}-a_{1})-
(x_{0}-a_{2}),...,(x_{0}-a_{1})-(x_{0}-a_{n})\}=\bR^{d},
\end{equation}
 any $\xi$ has the said property, and hence $\lambda=\mu$.

Thus, for each point $x_{0}\in P^{o}$ there is only one
supporting plane at $(x_{0},U(x_{0}))$. This and the concavity
of $U$ implies that $U$ is continuously differentiable,
$\lambda=\nabla U(x_{0})$, and the lemma is proved.

\begin{corollary}
                              \label{corollary 5.17.1}

Take any representation
$$
x=\sum_{k}q_{k}a_{k}\quad\text{with}\quad\sum_{k}q_{k}=1.
$$
Then
$$
\sum_{k}\frac{q_{k}}{p_{k}(x)}=n.
$$
\end{corollary}

Indeed, it suffices to observe that   $\xi=0$
in \eqref{8.23.1}.

\begin{lemma}
The functions $U$, $p_{k}$ are infinitely differentiable
in $P^{o}$.

\end{lemma}

Proof. Denote $\lambda(x)=\nabla U(x)$. Then
$$
p_{k}(x)=\frac{1}{n-(x-a_{k},\lambda(x))},
$$
and $\lambda(x)$ satisfies
$$
F(\lambda(x),x)=0,
$$
where
$$
F(\lambda,x)=(F^{i}(\lambda,x),i=1,...,d),\quad
F^{i}(\lambda,x)=\sum_{k}\frac{1}{n-(x-a_{k},\lambda
)}(x^{i}-a_{k}^{i}).
$$
We have
$$
\frac{\partial}{\partial\lambda^{j}}F^{i}(\lambda,x)
=\sum_{k}\frac{1}{(n-(x-a_{k},\lambda))^{2}}(x^{j}-a_{k}^{j})
(x^{i}-a_{k}^{i}).
$$
By \eqref{5.21.1} there is no nonzero vectors $\eta$ that are
orthogonal to all $x-a_{k}$. It follows that
 the matrix with the  entries
$\tfrac{\partial}{\partial\lambda^{j}}F^{i}(\lambda,x)$
 is nondegenerate,
$\lambda(x)$ is infinitely differentiable by  the implicit function
theorem and the lemma is proved.

\begin{lemma}

                                         \label{theorem 5.15.1}

Let $x\in P^{o}$, $\xi\in\bR^{d}$. Then
\begin{equation}
                                                    \label{11.15.1}
 U_{(\xi)(\xi)}(x)=-\sum_{k}\frac{(p_{k(\xi)}(x))^{2}}
{p_{k}^{2}(x)}.
\end{equation}
Furthermore, if $\xi=\sum_{k}q_{k}( a_{k}-x)$
for some numbers $q_{k}$, then
\begin{equation}
                                                    \label{11.15.2}
 -U_{(\xi)(\xi)}(x)
= \sum_{k}\frac{ p_{k(\xi)}(x)q_{k}}
{p_{k}^{2}(x)} -U_{(\xi)}(x)\sum_{k}q_{k},
\end{equation}
\begin{equation}
                                                    \label{11.16.1}
\sum_{k}\frac{(p_{k(\xi)}(x))^{2}}
{p_{k}^{2}(x)}
\leq (n+1)\sum_{k}\frac{ q_{k}^{2}}
{p_{k}^{2}(x)}.
\end{equation}

Finally, $| p_{k( a_{k}-x)}(x)|\leq  (n+1)^{1/2}$ for any $k
=1,...,n$.
\end{lemma}

Proof. By differentiating \eqref{5.31.2} we find
$$
\frac{p_{k(\xi)}(x) }
{p_{k}^{2}(x)}=-U_{(a_{k}-x)(\xi)}(x)+U_{(\xi)}(x).
$$
By multiplying this equality by $q_{k}$ and summing up with respect
to $k$ we get   \eqref{11.15.2}
provided that $\xi=\sum_{k}q_{k}( a_{k}-x)$.  
Differentiating $\sum_{k}p_{k}(x)a_{k}=x$ and $\sum_{k}p_{k}(x)=1$ 
yields
$$
\sum_{k}p_{k (\xi)}(x)a_{k}=\xi,\quad
\sum_{k}p_{k (\xi)}(x) =0,\quad
\xi=\sum_{k}p_{k(\xi)}(x)( a_{k}-x),
$$
 which allows us to use  \eqref{11.15.2}   with $q_{k}=p_{k(\xi)}$ 
and obtain \eqref{11.15.1}.

Next, the right-hand side of \eqref{11.15.2} equals
$$
\sum_{k}\frac{q_{k}}{p_{k}}\big(\frac{p_{k(\xi)}}{p_{k}}
-p_{k}U_{(\xi)}\big).
$$
Its square
by H\"older's inequality is less than
$$
\sum_{k}\frac{ q_{k}^{2}}
{p_{k}^{2} }\sum_{k}\big(\frac{ (p_{k(\xi)})^{2}}
{p_{k}^{2}}-2p_{k(\xi)}U_{(\xi)}+p_{k}^{2}U_{(\xi)}^{2}\big).
$$
We recall \eqref{11.15.1} and observe that
$$
\sum_{k}p_{k(\xi)}=0,\quad\sum_{k}p_{k}^{2}\leq1,
\quad
 U_{(\xi)}^{2}=\big(\sum_{k}\frac{p_{k(\xi)}}{p_{k}}\big)^{2}
\leq -n U_{(\xi)(\xi)} .
$$
Then we find that
$$
 U_{(\xi)(\xi)}^{2}\leq(n+1)|U_{(\xi)(\xi)}|\sum_{k}\frac{ q_{k}^{2}}
{p_{k}^{2} },
$$
which is equivalent to \eqref{11.16.1}.

The last assertion of the lemma is obtained by taking $\xi=a_{k}-x$
in \eqref{11.16.1}.  
The lemma is proved.  
\begin{theorem}
                                  \label{lemma 07.5.8.1}

For $\xi\ne0$ in  $  P^{o}$ we have
$$
 U_{(\xi)(\xi)}(x)\geq-(n+4n^{2})
\frac{|\xi|^{2}}{d^{2}(x,\xi)
\wedge d^{2}(x,-\xi)}.
$$
In particular, for any $x$ in $P^{o}$ we have
$$
\frac{|p_{k(\xi)}(x)|}{p_{k}(x)}\leq N\frac{|\xi|}{
d(x,\xi)
\wedge d(x,-\xi)},
$$
where $N=(n+4n^{2})^{1/2}$.
\end{theorem}

Proof. Without losing generality we 
assume that $|\xi|=1$, take $x\in
P^{o}$, and set
$$
y=x+d(x,\xi)\xi ,\quad y(t)=(1-t)x+ty,\quad t\in(-\varepsilon ,1),
$$
where $\varepsilon>0$ is to be chosen later.
Certainly there is a representation
$$
y =\sum_{k}q_{k}a_{k},\quad q_{k}\geq0,\quad
\sum_{k}q_{k}=1.
$$
Therefore, for sufficiently small $\varepsilon$
and all $t\in(-\varepsilon ,1)$ we have
$$
y(t)=\sum_{k}a_{k}((1-t)p_{k}(x)+tq_{k}),\quad 
(1-t)p_{k}(x)+tq_{k}>0,
$$
$$
\sum_{k}((1-t)p_{k}(x)+tq_{k})=1.
$$
By definition,
$$
U(y(t))\geq\sum_{k}\ln((1-t)p_{k}(x)+tq_{k})
$$
with equality for $t=0$. Therefore, the second
derivatives in $t$
at $t=0$ of the extreme terms are linked by a similar inequality,
 that is
$$
d^{2}(x,\xi)U_{(\xi)(\xi)}(x)\geq-\sum_{k}\frac{(p_{k}(x)-q_{k})^{2}}
{p_{k}^{2}(x)} \geq -n-\sum_{k}\frac{ q_{k} ^{2}}
{p_{k}^{2}(x)}.
$$
 
In like manner for $z=x-d(x,-\xi)\xi$ we find
$$
z=\sum_{k}r_{k}a_{k},\quad r_{k}\geq0,\quad\sum_{k}r_{k}=1,
$$
$$
d^{2}(x,-\xi)U_{(\xi)(\xi)}(x)\geq  -n-\sum_{k}\frac{ r_{k} ^{2}}
{p_{k}^{2}(x)}.
$$

However, for some $\alpha,\beta>0$ such that $\alpha y+\beta z=x$
and $\alpha+\beta=1$
we have
$$
x=\sum_{k}(\alpha q_{k}+\beta r_{k})a_{k},\quad
\sum_{k}(\alpha q_{k}+\beta r_{k})=1.
$$

 By Corollary \ref{corollary 5.17.1}
$$
\sum_{k}\frac{\alpha q_{k}+\beta r_{k}}{p_{k}(x)}=n.
$$
It follows that  
$$
\sum_{k}\frac{q_{k}^{2}}{p_{k}^{2}(x)}
\leq\big(\sum_{k}\frac{q_{k} }{p_{k}(x)}\big)^{2}
\leq \alpha^{-2}n^{2},\quad
\sum_{k}\frac{r_{k}^{2}}{p_{k}^{2}(x)}
\leq \beta^{-2}n^{2}
$$
and hence at least one of 
$$
\sum_{k}\frac{q_{k}^{2}}{p_{k}^{2}(x)},\quad
\sum_{k}\frac{r_{k}^{2}}{p_{k}^{2}(x)}
$$
is less than $4n^{2}$. This yields the result
and the theorem is proved. 

Now we are going to get prepared to estimating the Lipschitz
constants of $p_{k}$'s. Recall that $a_{1}=0$ and let
 $P_{1}$ be the polyhedron
with vertices $a_{2},..., a_{n}$, let $U_{n-1}(x)$ be the function
$U$ defined relative to $P_{1}$, and let $P_{1}^{o}$ be the relative
interior of $P_{1}$.

\begin{lemma}
                                       \label{lemma 5.20.1}
Let $x\in P^{o}$ and let $\lambda>1$
be such that $\lambda x\in P_{1}^{o}$. Then

(i)
\begin{equation}
                                                \label{5.31.1}
U_{n-1}(\lambda x) -n\ln\lambda+\ln(\lambda-1)
\leq U(x);
\end{equation}

(ii) we have an equality in \eqref{5.31.1} instead of the
inequality if we take $\lambda=\lambda(x):=(1-p_{1}(x))^{-1}$.
\end{lemma}

Proof. Let $\bar{p}_{2},...,\bar{p}_{n}$ be the set
that achieves $U_{n-1}(\lambda x)$. Then
\begin{equation}
                                                \label{8.23.2}
\sum_{i=2}^{n}\bar{p}_{i}=1,\quad\sum_{i=2}^{n}
\bar{p}_{i}a_{i}=\lambda x.
\end{equation}
Therefore, for $p_{i}:=\lambda^{-1}\bar{p}_{i}$, $i=2,...,n$,
and $p_{1}:=1-p_{2}-...-p_{n}=1-\lambda^{-1}$ we have
$p_{1}>0$, since $\lambda>1$, and
$$
\sum_{i=1}^{n}p_{i}=1,\quad\sum_{i=1}^{n}p_{i}a_{i}
=\sum_{i=2}^{n}\lambda^{-1}\bar{p}_{i}a_{i}=x.
$$
By adding that
\begin{equation}
                                                \label{8.23.3}
U_{n-1}(\lambda x)=\sum_{i=2}^{n}\ln\bar{p}_{i}
=(n-1)\ln\lambda+\sum_{i=1}^{n}\ln p_{i}-\ln(1-\lambda^{-1}).
\end{equation}
we certainly obtain \eqref{5.31.1}.

To prove assertion (ii) observe that for $\lambda=\lambda(x)$,
$p_{i}:=p_{i}(x)$,
and $\bar{p}_{i}:=\lambda p_{i}(x)$ we have
\eqref{8.23.2} and $\bar{p}_{i}>0$. It follows that
$\lambda x\in P_{1}^{o}$ and the first equality
sign in
\eqref{8.23.3} should be replaced with $\geq$. By combining this
with
\eqref{5.31.1} we get what we need. The lemma is proved.

\begin{corollary} 
                                       \label{corollary 5.18.1}

For $x\in P^{o}$ we have $\lambda(x)x\in
P_{1}^{o}$ and
the set $\lambda(x)p_{2}(x)$,..., $\lambda(x)p_{n}(x)$
achieves $U_{n-1}(\lambda(x)x)$, so that
if for $y\in P^{o}_{1}$ we denote by
$\bar{p}_{2}(y)$,..., $\bar{p}_{n}(y)$ the set that
achieves $U_{n-1}(y)$, then for $x\in P^{o}$ and $k\geq2$
we have
\begin{equation}
                                                \label{9.27.1}
p_{k}(x)=(1-p_{1}(x))\bar{p}_{k}(\frac{x}{1-p_{1}(x)}).
\end{equation}
\end{corollary}

\begin{theorem}
The functions $p_{1}(x),...,p_{n}(x)$
are Lipschitz continuous in $P^{o}$ and, therefore,
admit extensions to Lipschitz continuous functions
in  $P$.

\end{theorem}

Proof. We will be using the induction on $n$.  If $n=2$
and say $P=[0,1]$, $a_{1}=0,a_{2}=1$, then
$p_{1}(x)=1-x,p_{2}(x)=x$, and our assertion is true indeed.

Assume that our assertion is proved for all polyhedra
with $n-1$ vertices. Then the functions $\bar{p}_{k}$
introduced in Corollary \ref{corollary 5.18.1}
are Lipschitz continuous in $P_{1}^{o}$.
Since $|p_{1(x)}|\leq (n+1)^{1/2}$, from \eqref{9.27.1} we have
that for $k\geq2$ and $x\in P^{o}$
$$
|p_{k (x)}(x)|\leq N+(1-p_{1}(x))\lim_{\varepsilon\downarrow0}
\varepsilon^{-1}\big|\bar{p}_{k}(\frac{x+\varepsilon x }
{1-p_{1}(x+\varepsilon x)})-\bar{p}_{k}(\frac{x}{1-p_{1}(x)})\big|
$$
$$
\leq  N+N(1-p_{1}(x))
\lim_{\varepsilon\downarrow0}
\varepsilon^{-1}\big| \frac{x+\varepsilon x }
{1-p_{1}(x+\varepsilon x)} -\ \frac{x}{1-p_{1}(x)} \big|
$$
$$
=N+N|x|\big|1+\frac{p_{1(x)}(x)}{1-p_{1}(x)}\big|
\leq N+N\frac{|x|}{ 1-p_{1}(x) }.
$$

Thus, for any $\varepsilon>0$, $ p_{k(x)} (x)$ are bounded
as long as $|x|\geq\varepsilon$ and $k\geq2$. 
Above we also used  that $ p_{1(x)} (x)$ is bounded.

Generally,
$p_{k(x-a_{j})}(x)$ are bounded as long as $|x-a_{j}|\geq\varepsilon$.
In particular, $p_{1(x-a_{j})}(x)$ are bounded
for $|x-a_{j}|\geq\varepsilon$. We now claim that
there exists an $\varepsilon>0$ and $N_{0}$ such that, for any unit
$\xi\in\bR^{d}$ and $x\in P$ one can find 
numbers $\eta_{1},...,\eta_{k}$ such that
$$
\xi=\sum_{k:|x-a_{k}|\geq\varepsilon}\eta_{k}(x-a_{k}),
\quad|\eta_{k}|\leq N_{0}.
$$

Indeed, if $\varepsilon$ is small enough the restriction
of summation may exclude
only one term with $k$ such that $|x-a_{k}|<\varepsilon$.
Still the remaining set $\{x-a_{j},j\ne k\}$
would be close to $\{a_{k}-a_{j},j\ne k\}$ a subset of which
forms a basis in $\bR^{d}$. On the other hand, if there is nothing 
to exclude, our claim follows from \eqref{5.21.1}.

This proves that $|p_{1(\xi)}(x)|$ is bounded for $x\in P^{o}$,
$|\xi|=1$. Of course, the same holds for other $|p_{k(\xi)}(x)|$
and the theorem is proved.

\mysection{Proof of Theorem \protect\ref{theorem 07.5.7.1}}
                                            \label{section 5.21.3}

First we introduce a few new objects.
Let an integer $r\in [2,d]$ and let $\Gamma_{1},...,\Gamma_{r}\in\Phi$ be such that
$n_{\Gamma_{1}},...,n_{\Gamma_{r}}$ are linearly independent. Then
$$
|n_{\Gamma_{r}}-\Pi_{\text{\rm Span}\,(n_{\Gamma_{1}},...,n
_{\Gamma_{r-1}})}n_{\Gamma_{r}}|>0,
$$
where $\Pi_{\cL}$ is the orthogonal projection operator on
a subspace $\cL\in\bR^{d}$.
Since there are only finitely many such $r$ and 
$\Gamma_{1},...,\Gamma_{r}\in\Phi$, we see that there is 
a constant $\kappa\geq1$ such that we always have
$$
|n_{\Gamma_{r}}-\Pi_{\text{\rm Span}\,(n_{\Gamma_{1}},...,n
_{\Gamma_{r-1}})}n_{\Gamma_{r}}|\geq\kappa^{-1}.
$$

For a $\zeta_{1}>0$ define recursively
$$
\zeta_{r}=2\kappa\sum_{i=1}^{r-1}\zeta_{i}\quad r\geq2.
$$
Obviously, $\zeta_{r}\geq\zeta_{r-1}$ and $\zeta_{r}$ is  a linear function of
$\zeta_{1}$ so that we can choose and fix a $\zeta_{1}>0$ such that
$$
\gamma_{r}:=2^{r}\zeta_{r}\leq1/2,\quad r=1,...,d+1.
$$
Set 
$$
\varepsilon_{r}=\kappa^{-1}\gamma_{r}\quad (\leq\gamma_{r}).
$$

Now fix $x$ and drop it in some notation.
Take    a $\Gamma\in\Phi$ for
which
$d(x,\xi)\wedge d(x,-\xi)$ equals the distance  from $x$ to $\Gamma$ along the
line
$x+t\xi|\xi|^{-1}$, $t\in\bR$. 
Then
$$
d(x,\xi)\wedge d(x,-\xi)
=|\xi|\frac{d _{\Gamma}(x)}{|(n_{\Gamma},\xi)|}.
$$

Denote $\delta=|\nabla p_{k}|^{-1}$.
There are two cases.

{\em Case 1\/} 
$$
|(n_{\Gamma},\xi)|\geq \varepsilon_{1}
 \delta|p_{k(\xi)}| .
$$

{\em Case 2\/} 
$$
|(n_{\Gamma},\xi)|< \varepsilon_{1}
 \delta|p_{k(\xi)}| .
$$

In the first case additionally assume that
  $p_{k}(x)\leq d _{\Gamma}(x)$. Then
by Theorem~\ref{lemma 07.5.8.1}
$$
\frac{|p_{k(\xi)} |}{
p_{k}^{1/2} }\leq N|\xi|\frac{p_{k}^{1/2} }{
d(x,\xi)\wedge d(x,-\xi)}=N
\frac{p_{k}^{1/2} |(n_{\Gamma},\xi)|}{
d_{\Gamma} }\leq N
\frac{|(n_{\Gamma} ,\xi)|}{d^{1/2}_{\Gamma} }.
$$

On the other hand, if
$p_{k}(x)\geq d _{\Gamma}(x)$, then
$$
\frac{|p_{k(\xi)} |}{
p_{k}^{1/2} }\leq  (\varepsilon_{1}\delta)^{-1}
\frac{|(n_{\Gamma} ,\xi)|}{
p_{k}^{1/2} }
\leq (\varepsilon_{1}\delta)^{-1}
\frac{|(n_{\Gamma} ,\xi)|}{d^{1/2}_{\Gamma} }.
$$
Thus, in Case 1 we have
$$
\frac{|p_{k(\xi)} |}{
p_{k}^{1/2} }
\leq N(\varepsilon_{1}\delta)^{-1}
\frac{|(n_{\Gamma} ,\xi)|}{d^{1/2}_{\Gamma} },
$$
which proves \eqref{07.5.7.1} since $\delta^{-1}$
is a bounded function.

In the rest of the proof we concentrate on Case 2.
We will be using a recursive procedure.
Denote $\xi_{1}=\xi$, $\Gamma_{1}=\Gamma$, $n_{1}=n_{\Gamma_{1}}$,
and introduce 
$$
\xi_{2}=\xi_{1}-n_{1}(n_{1},\xi_{1}).
$$
Observe that (Case 2)
$$
|p_{k(\xi_{2})}-p_{k(\xi_{1})}|=
|p_{k(n_{1})}|\cdot|(n_{1},\xi_{1})|
\leq\delta^{-1}  |(n_{1},\xi_{1})|
<\varepsilon_{1}|p_{k(\xi_{1})}|
\leq\gamma_{1}|p_{k(\xi_{1})}|,
$$
$$
(1-\gamma_{1})|p_{k(\xi_{1})}|<
|p_{k(\xi_{2})}|<(1+\gamma_{1})|p_{k(\xi_{1})}|.
$$
In particular, $p_{k(\xi_{2})}\ne0$ and $\xi_{2}\ne0$.
Also in Case 2 we have $\xi_{2}\perp n_{1}$ and
$$
|\xi_{2}-\xi_{1}|=|(n_{1},\xi_{1})|\leq\varepsilon_{1}\delta|
p_{k(\xi_{1})}|\leq\gamma_{1}\delta|
p_{k(\xi_{1})}|.
$$

It follows that for an integer $r\geq2$ we have  
$\Gamma_{i}\in\Phi$, $i=1,...,r-1$, and vectors
$\xi_{i}\ne0$, $i=1,...,r $,  such that
for $n_{i}$ being the normal vectors to $\Gamma_{i}$
we have

(i) $n_{1},...,n_{r-1}$ are linearly independent;

(ii) for $1\leq j<i\leq r$ we have
$\xi_{i}\perp n_{j}$;

(iii) for $2\leq i\leq r$ we have 
$$
(1-\gamma_{i-1})|p_{k(\xi_{i-1})}|<
|p_{k(\xi_{i})}|<(1+
\gamma_{i-1})|p_{k(\xi_{i-1})}|;
$$

(iv) for $2\leq i\leq r$ we have
$|\xi_{i}-\xi_{i-1}|\leq\gamma_{i-1}
\delta|p_{k(\xi_{i-1})}|$;

(v) for $1\leq i\leq r-1$ we have
$|(n_{i},\xi_{i})|\leq \varepsilon_{i}
 \delta|p_{k(\xi_{i})}|$:

(vi) for $1\leq i\leq r-1$ the face $\Gamma_{i}\in\Phi$ is the one
for which $d(x,\xi_{i})\wedge d(x,-\xi_{i})$ equals the distance
 from $x$
to $\Gamma_{i}$ along the line
  $x+t\xi_{i}|\xi_{i}|^{-1}$, $t\in\bR$.

In light of (i) we have $r-1\leq d$. Also
observe that by virtue of (iii) and (iv) for $2\leq i\leq r$
(recall that $\gamma_{j}\leq1/2$)
\begin{equation}
                                                 \label{07.5.13.3}
|p_{k(\xi_{i-1})}|
\leq |p_{k(\xi_{r})}|\prod_{j=i-1}^{r-1}(1-\gamma_{j})^{-1}
\leq 2^{r-i+1}|p_{k(\xi_{r})}|,
\end{equation}
\begin{equation}
                                                 \label{07.5.13.2}
|\xi_{r}-\xi_{1}|\leq\sum_{i=2}^{r}
|\xi_{i}-\xi_{i-1}|\leq\delta
|p_{k(\xi_{r})}|\sum_{i=2}^{r}\gamma_{i-1}2^{r-i+1} 
=\varepsilon_{r}\delta|p_{k(\xi_{r})}| /2.
\end{equation}

Now introduce $\Gamma_{r}$
as the face of $P$
for which $d(x,\xi_{r})\wedge d(x,-\xi_{r})$ equals the distance
 from $x$
to $\Gamma_{r}$ along the line
  $x+t\xi_{r}|\xi_{r}|^{-1}$, $t\in\bR$. Set $n_{r}=n_{\Gamma_{r}}$
and first suppose that
\begin{equation}
                                                 \label{07.5.11.2}
|(n_{r},\xi_{r})|\geq\varepsilon_{r}
\delta|p_{k(\xi_{r})}|.
\end{equation}

Then as in Case 1 
\begin{equation}
                                                 \label{07.5.13.1}
\frac{|p_{k(\xi_{r})}|}{p_{k}^{1/2}}
\leq N
(\varepsilon_{r}\delta)^{-1}
\frac{|(n_{r},\xi_{r})|}{d^{1/2}_{\Gamma_{r}}}.
\end{equation}

Here by \eqref{07.5.13.3}
 the left-hand side dominates 
$$
2^{-r}\frac{|p_{k(\xi_{1})}|}{p_{k}^{1/2}}.
$$
To estimate the right-hand side of \eqref{07.5.13.1} use
\eqref{07.5.13.2} and \eqref{07.5.11.2} to get
$$
|(n_{r},\xi_{r})-(n_{r},\xi_{1})|
\leq|\xi_{r}-\xi_{1}|
\leq(1/2)|(n_{r},\xi_{r})|.
$$
 
Hence,
$$
|(n_{r},\xi_{r})|\leq2|(n_{r},\xi_{1})|
$$
and going back to \eqref{07.5.13.1} we obtain
$$
\frac{|p_{k(\xi_{1})}|}{p_{k}^{1/2}}
\leq N
\frac{|(n_{r},\xi_{1})|}{d^{1/2}_{\Gamma_{r}}},
$$
which proves \eqref{07.5.7.1}.

In the situation that \eqref{07.5.11.2} is violated
 introduce
$$
\xi_{r+1}=\xi_{r }-f_{r }(n_{r },\xi_{r }) ,
$$
where $f_{r }=h_{r }/|h_{r }|^{2}$ and
$$
h_{r }=n_{r }-\Pi_{\text{\rm Span}\,(n_{1},...,
n_{r-1})}n_{r }.
$$
Observe that $h_{r}\ne0$. Otherwise, $n_{r}$ would lie
in $\text{\rm Span}\,(n_{1},...,
n_{r-1})$, $\xi_{r}$ would be orthogonal also to $n_{r}$
and the line $x+t\xi_{r}|\xi_{r}|^{-1}$, $t\in\bR$, would have never met
$\Gamma_{r}$. In particular, property (i) holds
with $r+1$ in place of $r$. By the definition of $\kappa$
 we have $|h_{r}|\geq
\kappa^{-1}$. Then $|f_{r}|\leq\kappa$ and since
$|(n_{r},\xi_{r})|<\varepsilon_{r}\delta |p_{k(\xi_{r})}|$, we have
$$
|p_{k(\xi_{r+1})}-p_{k(\xi_{r})}|=|p_{k(f_{r})}|
\cdot|(n_{r},\xi_{r})|< \kappa
\varepsilon_{r}|p_{k(\xi_{r})}|=\gamma_{r}|p_{k(\xi_{r})}|
$$
implying that 
$p_{k(\xi_{r+1})}\ne0$, $\xi_{r+1}\ne0$ and (iii) holds with $r+1$ in
place of $r$. Also notice that, for $j\leq r-1$, we have
 $\xi_{r},h_{r},f_{r}\perp n_{j}$, which implies that $\xi_{r+1}\perp n_{j}$.
Furthermore,
$$
(h_{r},n_{r})=(h_{r},n_{r}-
\text{\rm Proj}_{\text{\rm Span}\,(n_{1},...,
n_{r-1})}n_{r })=|h_{r}|^{2},\quad (f_{r},n_{r})=1,
$$
$$
(\xi_{r+1},n_{r})=(\xi_{r},n_{r})-(f_{r},n_{r})
(n_{r},\xi_{r})=0,
$$
so that $\xi_{r+1}\perp n_{j}$ for all $j\leq r$
and (ii) holds with $r+1$ in place of $r$.
Properties (v) and (vi) hold with $r+1$ in place of $r$
by the assumption and construction.

Finally, 
$$
|\xi_{r+1}-\xi_{r}|=|f_{r}|\cdot|(n_{r},\xi_{r})|
\leq\kappa\varepsilon_{r}\delta|p_{k(\xi_{r})}|
=\gamma_{r}\delta|p_{k(\xi_{r})}|,
$$
so that (iv) holds with $r+1$ in place of $r$.

Thus, if \eqref{07.5.11.2}
 is violated, we can find objects $\Gamma_{i},\xi_{i}\ne0$
having the properties (i)-(vi) with $r+1$ in place of $r$.
This recursive process will stop at least when $r$ reaches $d+1$,
just because property (i) will prevent us from finding
$n_{r+1}$, which implies that at least at this moment \eqref{07.5.11.2}
should be satisfied. This proves the theorem.

\mysection{Additional information}
                                            \label{section 5.21.4}
\begin{remark}
                                            \label{remark 5.20.1}

One can estimate $p_{k}(x)$ from below for $x\in P^{o}$. It turns out
that
\begin{equation}
                                                    \label{5.20.1}
p_{k}(x)\geq\frac{d(x,x-a_{k})}{nd(x,x-a_{k})+n|x-a_{k}|},
\end{equation}
where $nd(x,x-a_{k})+n|x-a_{k}|$ 
is obviously bounded away from zero.
This and the fact that  $|\nabla p_{k}(x)|$
is bounded,
actually, show that $p_{k}(x)$ behaves like $ d(x,x-a_{k})$.

Indeed, take any $\xi\in\bR^{d}\setminus\{0\}$ and
observe that
$y:=x+\xi d(x,\xi)/|\xi|\in P$ can be written as
$$
y=\sum_{k}q_{k}a_{k},\quad q_{k}\geq0,\quad\sum_{k}q_{k}=1.
$$
Then
$$
\xi d(x,\xi)/|\xi|=y-x=\sum_{k}(p_{k}(x)-q_{k})(x-a_{k})
$$
and by \eqref{8.23.1}
$$
U_{(\xi)}(x) d(x,\xi)/|\xi|=\sum_{k}
\frac{q_{k}-p_{k}(x)}{p_{k}(x)}=
\sum_{k}
\frac{q_{k} }{p_{k}(x)}-n\geq -n,
$$
so that, for any $\xi\ne0$
$$
 U_{(\xi)}(x)\geq-\frac{n|\xi|}{d(x,\xi)}.
$$
If $\xi=x-a_{k}$, this  and \eqref{5.31.2}
imply that
$$
\frac{1}{p_{k}(x)}=n-U_{(x-a_{k})}(x)\leq n+
\frac{n|x-a_{k}|}{d(x,x-a_{k})},
$$
which is equivalent to \eqref{5.20.1}.

\end{remark}

\begin{remark}
If $\xi=\sum_{k}q_{k}( a_{k}-x)$ and $\sum_{k}q_{k}=0$ then
\begin{equation}
                                             \label{5.20.7}
\sum_{k}\frac{ q_{k}^{2}}
{p_{k}^{2}(x)}=
\sum_{k}\frac{(p_{k(\xi)}(x))^{2}}
{p_{k}^{2}(x)}
  +\sum_{k}\frac{(p_{k(\xi)}(x)-q_{k})^{2}}
{p_{k}^{2}(x)} .
\end{equation}

Indeed,  write
$$
\sum_{k}\frac{ q_{k}^{2}}
{p_{k}^{2}(x)}=\sum_{k}\frac{(q_{k}- p_{k(\xi)} )^{2}}
{p_{k}^{2}(x)}+2\sum_{k}\frac{p_{k(\xi)}(q_{k}- p_{k(\xi)} )}
{p_{k}^{2}(x)}+\sum_{k}\frac{ (p_{k(\xi)})^{2}}
{p_{k}^{2}(x)}
$$
and observe that the middle term on the right is zero
due to \eqref{11.15.1} and \eqref{11.15.2}.
\end{remark}

\begin{remark}
One can improve the estimate of $|p_{k( x-a_{k})}|$
from Lemma \ref{theorem 5.15.1}. It turns out that
\begin{equation}
                                                 \label{5.20.5}
1\geq p_{k( x-a_{k})}+1-p_{k}=\alpha_{k}p_{k}^{2},
\end{equation}
where
$$
\alpha_{k}= \frac{(p_{k(x-a_{k})}+1-p_{k})^{2}}{p_{k}^{2}}+
\sum_{i\ne k}\frac{(p_{i(x-a_{k})}-p_{i})^{2}}{p^{2}_{i}} .
$$
In particular, $p_{k}\geq p_{k(x-a_{k})}\geq p_{k}-1$.

Indeed,  we may concentrate on proving \eqref{5.20.5}
only for $k=1$ in which case we
 apply  
\eqref{5.20.7} with $\xi=x-a_{1}$.
Since 
$$
x-a_{1}=(p_{1}-1)(a_{1}-x)+p_{2}(a_{2}-x)+...+p_{n}(a_{n}-x)
$$
one can take $q_{1}=p_{1}-1,q_{2}=p_{2},...,q_{n}=p_{n}$. Then
by \eqref{11.15.1} and \eqref{5.20.7}
(recall that $a_{1}=0$)
$$
n-1+\frac{(p_{1}-1)^{2}}{p_{1}^{2}}=-U_{(x)(x)}
+\alpha_{1}.
$$

  On the other hand,  
by differentiating \eqref{5.31.2}
we find
$$
\frac{p_{1(x)}}{p_{1}^{2}}=U_{(x)(x)}+U_{(x)} 
=U_{(x)(x)}+n-\frac{1}{p_{1} }.
$$
Hence,
$$
n-1+\frac{(p_{1}-1)^{2}}{p_{1}^{2}}=-\frac{p_{1(x)}}{p_{1}^{2}}
+n-\frac{1}{p_{1} }+\alpha_{1},
$$
$$
-2p_{1}+1=-p_{1(x)}  -p_{1}+\alpha_{1} p_{1}^{2},
$$
and the equality in \eqref{5.20.5} follows. The   first inequality
follows from the fact that $\alpha_{1} p_{1}^{2}
\geq(p_{1(x)}+1-p_{1})^{2}$.

\end{remark}

\begin{remark} 
                                            \label{remark 5.20.2}
On can combine Remark \ref{remark 5.20.1}
and the fact that
  for any $j,k$
$$
|p_{k(x-a_{j})}(x)|\leq (n+1)^{1/2}p_{k}(x)/p_{j}(x).
$$
in order to investigate the behavior of $\nabla p_{i}(x)$
as $x$ approaches $\partial P$. The above mentioned fact follows from
\eqref{11.16.1} when $\xi=x-a_{j}$.
 \end{remark}
\begin{remark} The functions $p_{k}$ have a peculiar
symmetry. It turns out that for all $j,k$
$$
\frac{p_{k(x-a_{j})}}{p_{k}^{2}}+\frac{1}{p_{j}}
=
\frac{p_{j(x-a_{k})}}{p_{j}^{2}}+\frac{1}{p_{k}}.
$$

Indeed, differentiating \eqref{5.31.2} easily yields
\begin{equation}
                                   \label{9.27.2}
\frac{p_{k(x-a_{j})}}{p_{k}^{2}}+\frac{1}{p_{j}}-n
=U_{(x-a_{k})(x-a_{j})},
\end{equation}
where the right-hand side is symmetric with respect
to $j,k$. One can combine \eqref{9.27.2}
with Remarks \ref{remark 5.20.1} and \ref{remark 5.20.2}
to obtain some information about the behavior
of the second-order derivatives of $U$ near the boundary of $P$.
For instance, $|U_{(x-a_{k})(x-a_{j})}|\leq n
+2(n+1)^{1/2}(p_{k}p_{j})^{-1}$.
 \end{remark}

\begin{remark}  Lemma \ref{lemma 5.15.1} allows one to obtain
a precise information about the behavior of the 
first-order derivatives of
$U$ near the vertices of $P$. Indeed, \eqref{5.31.2} shows that
$$
\lim_{x\to a_{k}} U_{(x-a_{k})}(x)=n-1.
$$
\end{remark}

\begin{remark}
If in the situation of Lemma \ref{lemma 5.20.1}
we have $\lambda(x)x\in P^{o}$ then the derivative
with respect to $\lambda$ of the left-hand side of
\eqref{5.31.1} is zero at $\lambda=\lambda(x)$.
In this case, by substituting $\lambda=\lambda(x)$
into \eqref{5.31.1} and differentiating with respect to $x$
we find
$$
\lambda(x)U_{n-1,x^{i}}(\lambda(x)x)=U_{x^{i}}(x),\quad
i=1,...,d.
$$
In particular, the gradients of $U_{n-1}$ and $U$ are proportional
at corresponding points.
\end{remark}


\begin{thebibliography}{m}


\bibitem{BJ} G. Barles and E.R. Jakobsen, {\em Error bounds for
monotone approximation schemes for parabolic
Hamilton-Jacobi-Bellman equations\/}, Preprint.

\bibitem{BOZ}
J.F. Bonnans, E. Ottenwaelter, and H. Zidani, {\em A fast
algorithm for the two dimensional HJB equation 
of stochastic control\/},
M2AN Math. Model. Numer. Anal., Vol. 38 (2004), No. 4, 723-735.

\bibitem{DK}   Hongjie Dong and N.V. Krylov,
{\em  On the rate of convergence
of finite-difference approximations for parabolic Bellman equations
with Lipschitz coefficients in cylindrical
domains\/}, to appear in Applied Math. and Optimization.

\bibitem{Fr} M.I. Freidlin,   {\em The factorization of 
nonnegative definite matrices\/},
 Teor. Verojatnost. i Primenen.,  Vol. 13  (1968),
No. 2, 375-378 in Russian.

\bibitem{GT} D. Gilbarg and N.S. Trudinger,
``Elliptic partial differential equations of second order'',
2nd edition, Grundlehren der mathematischen
Wissenschaften, Vol. 224,
 Springer, Berlin-Heidelberg-New York-Tokyo, 1983.

\bibitem{Kr1}  N.V. Krylov,
  ``Nonlinear elliptic and parabolic
equations of second
  order'',   Reidel, Dordrecht, 1987.

\bibitem{Kr} N.V. Krylov, {\em A priori estimates
of smoothness of solutions to  difference  Bellman's equations 
  with linear and quasilinear operators\/}, 
Math. Comp., Vol. 76 (2007), 669-698. 

\bibitem{KD} H.J. Kushner  and
 P.G. Dupuis,  ``Numerical methods
for stochastic control problems
in continuous time", 2nd edition,
Springer Verlag, 2001.

\bibitem{MW} T. Motzkin and W. Wasow, {\em On the approximation of
linear
elliptic differential equations by difference
equations with positive coefficients\/}, J. Math and Phys., Vol. 31
(1952), 253-259.

\bibitem{PS} R. S. Phillips and L. Sarason, {\em Elliptic-parabolic 
equations of the second order\/},
 J. Math. Mech., Vol.  17 (1968), No. 2, 891-917.
\end{thebibliography}
\end{document}